\documentclass[10pt]{article}
\begin{document}      
\centerline {{\Large\bf Invariant and evolutionary properties}}  
\centerline {{\Large\bf of the skew-symmetric differential forms}}  
\centerline {L.I. Petrova}
\centerline{{\it Moscow State University, Russia, e-mail: ptr@cs.msu.su}}
\bigskip

The present work pursues the aim to draw attention to unique possibilities 
of the skew-symmetric differential forms. 

At present the theory of skew-symmetric {\it exterior} 
differential forms that possess invariant properties has been developed.

In the work the readers are introduced to the skew-symmetric 
differential forms that were called evolutionary ones because they possess 
evolutionary properties. 

The combined mathematical apparatus of exterior and  evolutionary 
skew-symmetric differential forms in essence is a new mathematical language. 
{\it This apparatus can describe transitions from nonconjugated operators to 
conjugated ones}. There are no such possibilities in any mathematical formalism. 

In the present work it has been shown that the properties of exterior or 
evolutionary forms are explicitly or implicitly accounted for in all 
mathematical (and physical) formalisms.

\bigskip
{\large\bf Introduction}

While studying the general equations of dynamics A.Poincare showed an existence 
of integral invariants. When investigating the integrability conditions of 
the differential equations and the sets of equations in total differentials 
E.Cartan discovered the invariant properties of integrands in multiple integrals 
and elucidated their own significance. He called them as the exterior differential 
forms because they obey the rules of the exterior multiplication by Grassmann 
(the skew-symmetry conditions). It appears that the exterior differential 
forms together with the operation of exterior differentiation may possess 
the invariant, group, tensor, structural and other properties that are of great 
functional and utilitarian importance. The exterior differential forms were 
found wide application in differential geometry and algebraic topology. 

The anasysis of differential equations shows that, besides of the 
skew-symmet\-rical differential forms, which possess invariant properties 
and have been named the exterior differential forms, there are other 
skew-symmetrical differential forms, which possess evolutionary properties. 
The author has named these forms as evolutionary differential forms.

The evolutionary forms, as well as the exterior forms, are differential 
forms with exterior multiplication.   
A radical  distinction between the evolutionary forms and the exterior 
ones consists in the fact that the exterior differential forms are defined on 
manifolds with {\it closed metric forms}, whereas the evolutionary differential 
forms are defined on manifolds with {\it unclosed metric forms}. 

This leads to the fact that the evolutionary forms and exterior ones possess 
the opposing properties and the opposing mathematical apparatus. Whereas the 
mathematical apparatus of exterior differential forms involves the identical 
relations and nondegenerate transformations, the mathematical apparatus of 
evolutionary forms involves  the nonidentical relations and degenerate 
transformations. Thus, they complement each other and make up some unified 
whole. Between the exterior and evolutionary forms there exists a connection. 
{\it The evolutionary differential forms generate the exterior differential 
forms}.

Owing to these properties the mathematical apparatus of exterior and 
evolutionary differential forms allows  description of discrete transitions, 
quantum steps, evolutionary processes,  generation of various structures.
 
A transition from evolutionary forms to  closed exterior forms 
describes a transition from nonconjugated operators to conjugated ones.
 
\section{Properties and specific features of the exterior differential 
forms}

It should be noted that this work is not intended to present the general 
theory of exterior differential forms. The reader can find information on 
this subject, for example, in [1-8]. Some initial information concerning 
exterior differential forms is outlined, the basic properties of the 
exterior differential forms and peculiarities of their mathematical apparatus 
are described. It is shown that the invariant properties of closed exterior 
forms reveal themselves in many branches of mathematics.

\bigskip
The exterior differential form of degree $p$ ($p$-form on the differentiable 
manifold) can be written as [5,7,8]
$$
\theta^p=\sum_{i_1\dots i_p}a_{i_1\dots i_p}dx^{i_1}\wedge
dx^{i_2}\wedge\dots \wedge dx^{i_p}\quad 0\leq p\leq n\eqno(1.1) 
$$
Here $a_{i_1\dots i_p}$ are the functions of the variables $x^{i_1}$,
$x^{i_2}$, \dots, $x^{i_p}$, $n$ is the dimension of space,
$\wedge$ is the operator of exterior multiplication, $dx^i$,
$dx^{i}\wedge dx^{j}$, $dx^{i}\wedge dx^{j}\wedge dx^{k}$, \dots\
is the local basis which satisfies the condition of exterior
multiplication:
$$
\begin{array}{l}
dx^{i}\wedge dx^{i}=0\\
dx^{i}\wedge dx^{j}=-dx^{j}\wedge dx^{i}\quad i\ne j
\end{array}\eqno(1.2)
$$
[From here on the symbol $\sum$ will be omitted and it will be 
implied that a summation is performed over double indices.  Besides, the 
symbol of exterior multiplication will be also omitted for the  
sake of presentation convenience].

The differential of the (exterior) form $\theta^p$ is expressed as 
$$
d\theta^p=\sum_{i_1\dots i_p}da_{i_1\dots
i_p}dx^{i_1}dx^{i_2}\dots dx^{i_p} \eqno(1.3)
$$
and is the differential form of degree $(p+1)$. \{The exterior  
differentiating operator $d$ allows one to pass on from the $p$-fold 
integral over the $p$-dimension closed manifold to the $(p+1)$-fold 
integral over the $(p+1)$-dimension manifold bounded by first one 
[1]\}.

Local domains of manifold are the basis of the exterior 
form. In this section we will consider the domains of the 
Euclidean space [9] or differentiable manifolds [8]. (Manifolds, on which 
the exterior differential forms and the evolutionary forms may be 
defined, and the influence of the manifold properties on the 
differential forms  will be discussed in more detail in section 2).

Let us consider some examples of the exterior differential form whose 
basis are the Euclidean space domains.

We consider the 3-dimensional space. In this case the differential 
forms of zero-, first- and second degree can be written as [7]  
$$\theta^0=a,\eqno(1.4)$$
$$\theta^1=a_1 dx^1+a_2 dx^2+a_3 dx^3,\eqno(1.5)$$
$$\theta^2=a_{12}dx^1 dx^2+a_{23}dx^2 dx^3+a_{31}dx^3 dx^1\eqno(1.6)$$
With account for conditions (1.2) their differentials are the forms
$$d\theta^0=\frac{\partial a}{\partial x^1}dx^1+\frac{\partial a}{\partial
x^2}dx^2+\frac{\partial a}{\partial x^3}dx^3,\eqno(1.7)$$
$$d\theta^1=\left(\frac{\partial a_2}{\partial x^1}-\frac{\partial
a_1}{\partial x^2}\right)dx^1 dx^2+\left(\frac{\partial a_3}{\partial
x^2}-\frac{\partial a_2}{\partial x^3}\right)dx^2 dx^3+$$
$$\left(\frac{\partial a_1}{\partial x^3}-\frac{\partial a_3}{\partial
x^1}\right)dx^3 dx^1, \eqno(1.8)$$
$$d\theta^2=\left(\frac{\partial a_{23}}{\partial
x^1}+\frac{\partial a_{31}}{\partial x^2}+\frac{\partial
a_{12}}{\partial x^3}\right)dx^1 dx^2 dx^3 \eqno(1.9)$$

From (1.4)--(1.9) one can see the following.

a) Any function is the form of zero degree. Its basis is a surface of zero 
dimension, namely, a variety of points. A differential 
of the form of zero degree is an ordinary differential of the 
function.

b) The form of first degree is a differential expression. An ordinary 
differential of the function is an example of the first-degree form. 

c) Coefficients of the differentials of the forms of the zero-, first- and 
second degrees give the gradient, curl, and divergence, 
respectively. That is, the operator $d$, referred to as the 
exterior differentiation, is some abstract generalization of the 
ordinary operators of gradient, curl,  and divergence. At this point 
it should be emphasized that in mathematical analysis the ordinary 
concepts of gradient, curl, and divergence are the operators 
applied to the vector, and in the theory of exterior forms 
gradient, curl, and divergence obtained as the results of exterior 
differentiating (the forms of the zero- first- and second degrees) are 
the operators applied to the pseudovector (the axial vector).

These examples are  presented in order to give some initial illustrative 
presentation of the exterior differential forms.

From these examples one can assure oneself that, firstly, the 
differential of the exterior form is also the exterior form 
(but with the degree greater by one), and, secondly, he can see 
that the components of the differential form commutator are 
coefficients of the form differential. Thus, a differential 
of the first-degree form $\omega=a_i dx^i$ can be written as 
$d\omega=K_{ij}dx^i dx^j$ where $K_{ij}$ are  components of the 
commutator for the form $\omega$ that are defined as 
$K_{ij}=(\partial a_j/\partial x^i-\partial a_i/\partial x^j)$.

\subsection*{Closed exterior differential forms} 

Closed differential forms possess the invariant properties.
 
A form is called a closed one if its differential is equal to zero:
$$
d\theta^p=0\eqno(1.10)
$$

From condition (1.10) one can see that the closed form \{the 
kernel of the operator $d$\} is a conservative 
quantity. (This means that it can correspond to the conservation law, namely, 
to some conservative quantity).

A differential of the form is a closed form. That is, 
$$
dd\omega=0\eqno(1.11)
$$
where $\omega$ is an arbitrary exterior form.

The form that is a differential of some other form \{a mapping of 
the operator $d$\}: 
$$
\theta^p=d\theta^{p-1}\eqno(1.12)
$$
is called an {\it exact} form. Exact forms prove to be closed 
automatically [4]
$$
d\theta^p=dd\theta^{p-1}=0\eqno(1.13)
$$

Here it is necessary to pay attention to the following points. In the 
above presented formulas it was implicitly assumed that the differential 
operator $d$ is the total one (that is, the 
operator $d$ acts everywhere in the vicinity of the point 
considered locally),  and therefore, it acts on the manifold of the 
initial dimension $n$. However, a differential may be internal. Such a 
differential acts on some structure with the dimension being less 
than that of the initial manifold. The structure, on which the exterior 
differential form may become a closed {\it inexact} form, is a pseudostructure 
with respect to its metric properties. \{Cohomology, sections of cotangent 
bundles, the eikonal surfaces, the characteristical and potential surfaces, 
and so on may be regarded as examples of pseudostructures.\} 
The properties of pseudostructures will be considered later.

If the form is closed on  pseudostructure only, the closure 
condition is written as
$$
d_\pi\theta^p=0\eqno(1.14)
$$
And the pseudostructure $\pi$ is defined from the condition 
$$
d_\pi{}^*\theta^p=0\eqno(1.15)
$$
where ${}^*\theta^p$ is the dual form. 
(For the properties of dual forms see [10]).

\{A skew-symmetric tensor corresponds to the exterior 
differential form on the differentiable manifold, and the pseudotensor that 
is dual to the skew-symmetric tensor corresponds to the dual form\}.

From conditions (1.14) and (1.15) one can see that the form 
closed on pseudostructure (a closed inexact form) is a conservative object, 
namely, this quantity conserves on pseudostructure. (This can also correspond 
to some conservation law, i.e. to conservative object).

The exact form is, by definition, a differential (see condition (1.12)). 
In this case the differential is total. The closed inexact form is 
a differential too. The closed inexact form is an interior (on 
pseudostructure) differential, that is
$$
\theta^p_\pi=d_\pi\theta^{p-1}\eqno(1.16)
$$

And so, any closed form is a differential of the form of a lower 
degree: the total one $\theta^p=d\theta^{p-1}$ if the form is exact, 
or the interior one $\theta^p=d_\pi\theta^{p-1}$ on pseudostructure if 
the form is inexact. (This may have the physical meaning: the form of lower 
degree may correspond to the potential, and the closed form by itself 
may correspond to the potential force.) 

From conditions (1.12) and  (1.16) it is possible to see that 
a relation between closed forms of different 
degree can exist.

Similarly to the differential relation between the exterior 
forms of sequential degrees, there is an integral connection. The relevant 
integral relation has the form [10]
$$
\int\limits_{c^{p+1}}d\theta^p=\int\limits_{\partial
c^{p+1}}\theta^p\eqno(1.17)
$$
In particular, the integral theorems by Stokes and Gauss follow 
from the integral relation for $p=1,2$ in three-dimensional space. 
\{From this relation one can see that an integral of the closed form 
over the closed curve vanishes (in the case of a smooth manifold). 
However, in the case of a complex 
manifold (for example, a not simply connected manifold with the homology 
class being nonzero) the integral of the closed form (in this case the 
form is inexact) over the closed curve is nonzero. It may be equal to a 
scalar multiplied by $2\pi$, which in this case corresponds, for 
example, to such a physical quantity as the charge [10]. Just such integrals 
are considered in the theory of residues.\}

\subsection*{Invariant properties of closed exterior differential 
forms. Conjugacy of exterior differential forms}

Since the closed form is a differential (a total one if the form is exact, 
or an interior one on the pseudostructure if the form is inexact), then it is 
obvious that the closed form proves to be invariant under all 
transformations that conserve a differential. 

The examples of such 
transformations are unitary, tangent, canonical, 
gradient transformations.  

A closure of exterior differential forms and hence their 
invariance result from the conjugacy of elements of exterior or dual forms.

From definition of the exterior differential form one can see that the 
exterior differential forms have complex structure. Specific 
features of the exterior form structure are homogeneity with 
respect to the basis, skew-symmetry, integrating  terms 
each including two objects of different nature 
(the algebraic nature for form coefficients, and the geometric nature 
for  base components). Besides,  the exterior form depends 
on the space dimension and on the manifold topology. The closure 
property of the exterior form means that any objects, namely, 
elements of the exterior form, components of elements, elements of 
the form differential, exterior and dual forms and others, turn 
out to be conjugated. It is a conjugacy that leads to realization of 
invariant and covariant properties of the exterior and dual 
forms that have a great functional and utilitary importance. A  
variety of objects of conjugacy leads to the fact that the closed forms 
can describe a great number of different 
structures, and this fact once again emphasizes great mathematical 
potentialities of the exterior differential forms.

Let us consider some types of conjugacy that make the exterior 
differential and dual formsn closed, that is, 
they make the form differentials equal to zero.

As it was pointed out before, components of the exterior 
form commutator are  coefficients of differential of this form. If the 
commutator of the form vanishes, the form differential vanishes 
too, and this indicates that the form is a closed one. Therefore, 
the closure property of the form may be recognized by finding 
whether or not the commutator of the form vanishes.

One of the types of conjugacy is that for the form coefficients.

Let us consider the exterior differential form of the first degree 
$\omega=a_i dx^i$. In this case,  
as it was pointed before, 
the differential will be expressed 
as $d\omega=K_{ij}dx^i dx^j$, where
$K_{ij}=(\partial a_j/\partial x^i-\partial a_i/\partial x^j)$ are 
components of the form commutator.

It is evident that the differential may vanish if the components 
of commutator vanish. One can see that 
the components of commutator $K_{ij}$ may vanish if derivatives of 
the form coefficients vanish. This is a trivial 
case. Besides, the components $K_{ij}$ may vanish if the 
coefficients $a_i$ are derivatives of some function $f(x^i)$, 
that is, $a_i=\partial f/\partial x^i$. In this case, 
components of the commutator are equal to a difference of the mixed 
derivatives
$$
K_{ij}=\left(\frac{\partial^2 f}{\partial x^j\partial
x^i}-\frac{\partial^2 f}{\partial x^i\partial x^j}\right)
$$
and, therefore, they vanish. One can see that the form 
coefficients $a_i$, which satisfy these conditions,  are 
conjugated quantities (the operators of mixed differentiation turn out 
to be commutative).

Let us consider the case when the exterior form is written as
$$
\theta=\frac{\partial f}{\partial x}dx+\frac{\partial f}{\partial y}dy
$$
where $f$ is the function of two independent variables $(x,y)$. It is evident 
that this form is closed because it is 
equal to the differential $df$. And, for the dual form
$$
{}^*\theta=-\frac{\partial f}{\partial y}dx+\frac{\partial
f}{\partial x}dy
$$
to be closed also,  it is necessary that its commutator be equal to zero
$$
\frac{\partial^2 f}{\partial x^2}+\frac{\partial^2 f}{\partial
y^2}\equiv \Delta f=0
$$
where $\Delta$ is the Laplace operator. As a result the function $f$ has 
to be harmonic one.

Assume that the exterior differential form of the first degree has 
the form $\theta=udx+vdy$, where  $u$ and $v$ are functions of  
two variables $(x,y)$. In this case, the closure condition of the 
form, that is, the condition under which the form commutator vanishes, 
takes the form
$$
K=\left(\frac{\partial v}{\partial x}-\frac{\partial u}{\partial
y}\right)=0
$$
One can see that this is one of the Cauchy-Riemann conditions for 
complex functions. The closure condition of the relevant dual 
form ${}^*\theta=-vdx+udy$ is the second Cauchy-Riemann condition. 
\{Here one can see a connection with functions of complex 
variables. If we consider the function $w=u+iv$ of the complex variables 
$z=x+iy$, $\overline{z}=x-iy$ that obeys the Cauchy-Riemann 
conditions, then to this function will correspond the closed exterior and dual 
forms. (The Cauchy-Riemann conditions are the 
conditions under which the function of complex variables does not 
depend on the conjugated coordinate $\overline{z}$). And to the 
harmonic function of complex variables there corresponds the closed 
exterior differential form, whose coefficients $u$ and $v$ are  
conjugated harmonic functions\}.

It can exist the conjugacy which makes the interior 
differential on the pseudostructure equal to zero, $d_\pi\theta=0$. Assume 
that the interior differential is the form of the first degree (the 
form itself is that of zero degree) and can be presented as 
$d_\pi\theta=p_x dx+p_y dy$, where $p$ is the form of zero 
degree (any function). Then the closure condition of the form  
($d_\pi\theta=p_x dx+p_y dy=0$) is
$$
\frac{dx}{dy}=-\frac{p_y}{p_x}\eqno(1.18)
$$
This is a conjugacy of the basis and derivatives of the form 
coefficients. One can see that formula (1.18) is one of the 
formulas of canonical relations [11]. The second formula of 
canonical relations follows from a condition that the dual form 
differential vanishes. This type of conjugacy  is 
connected with the canonical transformation. For a differential of the 
first degree form (in this case the differential is the form of 
the second degree) a corresponding transformation is a gradient one.

At this point it should be remarked that relation (1.18) is 
a condition that the implicit function to exist. That is, the 
closed (inexact) form of zero degree is the implicit function.  
\{This is an example  of connection between an exterior form and 
analysis\}.

A property of the exterior differential form being closed on 
pseudostructure points to another type of conjugacy, namely, a  
conjugacy of exterior and dual forms.

From the properties of closed forms analyzed above it becomes 
evident that there are exterior and interior types of conjugacy. 
To the interior type of conjugacy there are assigned a conjugacy 
of the form coefficients, coefficients of the form and the 
basis, conjugacy of  coefficients of the differential, and so on. 
To the exterior types of conjugacy there is assigned a 
conjugacy of the exterior and dual forms.  Another example of 
the exterior conjugacy is that between the forms of sequential  
degrees (the closed differential form is a differential of the  
form of less by one degree, that is, the closed form is conjugated 
to the form of less by one degree).

With the conjugacy another characteristic property of the exterior differential 
forms, namely, their duality, is connected that has the fundamental meaning. 
\{The conjugacy is some identical relation between 
two operators or mathematical objects.  A duality is a 
concept that means that one object carries a double meaning or 
that two objects with different meanings (of different physical 
nature) are identically connected. If one knows any dual object,  
he can obtain the other object\}. A conjugacy of objects of 
the exterior differential form generates a duality of the exterior forms.

The connection of the exterior and dual forms is an example of the 
duality. The exterior form and the dual form correspond to the 
objects of different nature: the exterior form corresponds to the 
physical (i.e. algebraic) quantity, and the dual form corresponds to some 
spatial (or pseudospatial) structure. At the same time, under 
conjugacy the duality of these objects manifests itself, that is, 
if one form is known, it is possible to find the other form. 
\{It will be shown below that the 
duality between exterior and dual forms elucidates a 
connection between physical quantities and spatial structures (and together 
they form a physical structure\}. 
Here the duality is also revealed in the fact that, if the degree of the 
exterior form equals $p$, the dimension of the spatial structure  
equals $N-p$, where $N$ is the space dimension.

Since the closed exterior form possesses invariant properties and the dual 
form corresponding to that possesses covariant properties, the 
invariance of the closed exterior form and the relevant covariance of 
the dual form is an example of the duality of exterior 
differential forms.

The other example of a duality of the closed form is connected 
with the fact that the closed form of degree $p$ is a 
differential of the form of degree $p-1$. This duality is manifested in 
that, on the one hand, as it was pointed out before, the closed exterior 
form is a conservative quantity, and on the other hand, the closed 
form can correspond to a potential force [12,13]. 
(In the works [12,13] the physical 
meaning of this duality has been illustarted, and it has been shown  
in respect 
to what the closed form manifests itself as a potential 
force and with what the conservative physical quantity is connected).

Here it should be emphasized that the duality is a property of the 
closed forms only. In the case of unclosed forms one cannot speak 
about duality.

The further manifestation of duality of exterior differential 
forms that needs for more attention is a duality of the concepts 
of closure and of integrability of differential forms.

The form that can be presented as a differential can be called 
an integrable one because it is possible to integrate it directly. 
In this context the closed form, which is a differential (total,  
if the form is exact, or interior, if the form is inexact), is 
integrable. And the closed exact form is integrable identically,  
whereas the closed inexact form is integrable on pseudostructure only.

The concepts of closure and integrability are not identical ones. 
The closure of the form is defined with respect to the form of 
greater by one degree, whereas the integrability is defined with respect 
to the form of degree less by one. Really, the form is closed one when 
the form differential, which is the form of greater by one degree, 
equals zero. And the form referred to as integrable one is a 
differential of some form of degree less by one. The closure and the 
integrability are dual concepts. Namely, the closure and the 
integrability are further examples of a duality of exterior 
differential forms.

Here it should emphasize once again that the duality is a 
property of closed forms only. In particular, nonclosure and 
nonintegrability are not dual concepts, and this will be revealed 
while analyzing evolutionary differential forms.

The duality of closed differential forms reveals under an availability 
of one or other type of conjugacy. The duality is a tool that untangles the 
mutual connection, the mutual changeability and the transitions 
between different objects. 

\subsection*{Differential and geometrical structure}

From the definition of a closed inexact exterior form one can see 
that to this form there correspond two conditions:

1) condition (1.14) is a closure condition of the exterior form itself, and

2) condition (1.15) is that of the dual form.

Conditions (1.14) and (1.15) can be regarded as equations for a binary object 
that combines the pseudostructure and the exterior differential form defined 
on this pseudostructure. Such a binary object can be named  Bi-Structure. 
The Bi-Structure combines both an algebraic object, namely, the closed 
exterior form, and a geometric one, namely, the pseudostructure as well.

In its properties this differential and geometrical structure is a well-known 
G-Structure. 
Here a new term Bi-Structure was introduced to distinguish it from G-structures 
to which there correspond closed inexact exterior forms.

The specific feature of this structure consists in the fact that it combines 
objects that possess duality.  The closed exterior differential form and 
the closed dual form are such objects. The existence of one  object implies 
that the other one exists as well. It does not seem to make sense to combine 
mathematically these two binary objects. However, this combined structure 
constitutes a unified whole that carries a double meaning. 
(This statement can be understand by analyzing the role of such structures in 
physical applications. This has been shown in the authors' works [12,13].)  

From conditions (1.14) and (1.15) one can see that Bi-Stricture constitute 
a conservative object, namely, a quantity that is conservative on the 
pseudostructure. Hence, such an object (Bi-Structure) can correspond to some 
conservation law. (It is from such object that the physical fields 
and corresponding manifolds are formed.)

The properties and specific features of Bi-Structures will be presented in 
section 2. The evolutionary forms allow us to describe characteristics 
of these structures. 

\bigskip
Invariant properties of closed exterior differential forms lie at the basis 
of mathematical apparatus of exterior forms. Below it will be described 
some specific features of the mathematical apparatus of exterior forms. 

\subsection*{Operators of the theory of exterior differential 
forms}
In  differential calculus the derivatives are basic elements of 
the mathematical apparatus. By contrast, the 
differential is an element of mathematical apparatus of the 
theory of exterior differential forms. It enables one to analyze  
the conjugacy of derivatives in various directions, which extends 
potentialities of differential calculus. 

The operator of exterior differential $d$ (exterior 
differential) is an abstract generalization of ordinary mathematical 
operations of the gradient, curl, and divergence in the vector 
calculus [7]. If, in addition to the exterior differential, we 
introduce the following operators: 1) $\delta$ for transformations that 
convert the form of $(p+1)$ degree into the form of $p$ degree, 2) $\delta'$ 
for cotangent transformations, 3) $\Delta$ for the 
$d\delta-\delta d$ transformation, 4) $\Delta'$ for the $d\delta'-\delta'd$ 
transformations, then in terms of these operators that act on the exterior 
differential forms one can write down the operators in the field 
theory equations. The operator $\delta$ corresponds to Green's 
operator, $\delta'$ does to the canonical transformation operator, 
$\Delta$ does to the d'Alembert operator in 4-dimensional space, and 
$\Delta'$ corresponds to the Laplace operator [8,10]. It can be seen 
that the operators of the exterior differential form theory are 
connected with many operators of mathematical physics.

\subsection*{Identical relations of exterior differential forms} 
In the theory of exterior differential forms the closed forms that 
possess various types of conjugacy play a principal role. Since 
the conjugacy is a certain connection between two operators or 
mathematical objects, it is evident that, to express a conjugacy mathematically, 
it can be used relations. Just such relations constitute the basis of 
mathematical apparatus of the exterior differential forms. This is 
an identical relation.

Identical relations for exterior differential forms reflect the 
closure conditions of differential forms, namely, vanishing the form 
differential (see formulas (1.10), (1.14), (1.15)) and the conditions 
connecting the forms of consequent degrees (see formulas (1.12), (1.16)).

An importance of the identical relations for exterior differential forms 
is manifested by the fact that practically in all branches of physics, 
mechanics, thermodynamics one faces such identical relations. 

One can present the following examples:  

a) the Poincare invariant [5] $ds\,=\,-H\,dt\,+\,p_j\,dq_j$,

b) the second principle of thermodynamics [14] $dS\,=\,(dE+p\,dV)/T$,

c) the vital force theorem in theoretical mechanics: $dT=X_idx^i$ 
where $X_i$ are components of  potential force, and $T=mV^2/2$ is the 
vital force,

d) the conditions on characteristics [11] in the theory of differential 
equations, and so on.
 
A requirement that the function is an antiderivative (the integrand is a 
differential of a certain function) can be written in terms of such 
identical relation.

An  existence of the harmonic function is written by means of the 
identical relation: the harmonic function is a closed 
form, that is a differential (a differential on the Riemann surface).

The identical relations in differential forms express the fact that each 
closed exterior form is a differential of some exterior form (with the degree 
less by one). In general form such an identical relation can be written as 
$$
d _{\pi}\phi=\theta _{\pi}^p\eqno(1.19)
$$
In this relation the form in the right-hand side has to be a {\it closed} 
one. (As it will be shown below, 
the identical relations are satisfied only on pseudostructures). 

In identical relation (1.19) in one side it stands the closed form and 
in other side does a differential of some differential 
form of the less by one degree, which is a closed form as well. 

In addition to relations in differential forms from  
the closure conditions of the differential forms and the conditions 
connecting the forms of consequent degrees  the identical relations of 
other types are obtained. The types of such relations are presented below. 

{\it Integral identical relations}.

The formulas by Newton, Leibnitz, Green, the integral relations by Stokes, 
Gauss-Ostrogradskii are examples of integral identical relations. 

{\it Tensor identical relations}.

From relations that connect exterior forms of consequent degrees 
one can obtain the vector and tensor identical relations that connect 
operators of the gradient, curl, divergence and so on.

From the closure conditions of exterior and dual forms one can obtain 
identical relations such as the gauge relations in electromagnetic field 
theory, the tensor relations between connectednesses and their derivatives 
in gravitation (the symmetry of connectednesses with respect to lower indices, 
the Bianchi identity, the conditions imposed on the Christoffel symbols) 
and so on.

{\it Identical relations between derivatives}.

The identical relations between derivatives correspond to the closure 
conditions of exterior and dual forms. The examples of such relations are 
the above presented Cauchi-Riemann conditions in the theory of complex 
variables, the transversality condition in the variational calculus, 
the canonical relations in the Hamilton formalism, 
the thermodynamic relations between derivatives 
of thermodynamic functions [14], a condition that a derivative of 
implicit function is subjected to, the eikonal relations [15] and so on. 

The above presented examples show that the identical relations for exterior 
differential forms occur in various branches of mathematics and physics.

Identical relations of exterior differential forms are the mathematical 
expression of various kinds of conjugacy that leads to closed exterior 
forms. They describe a conjugacy of any objects: the form elements, 
components of each element, exterior and dual forms, exterior forms of various 
degrees, and others. The identical relations, which are connected with different 
kinds of conjugacy, elucidate invariant, structural and group properties of 
exterior forms that are of great importance in applications. 

A functional significance of identical relations for exterior differential 
forms lies in the fact that they can describe a conjugacy of objects of 
different mathematical nature. This enables one to see internal connections 
between various branches of mathematics. Due to these possibilities the 
exterior differential forms have wide application in various branches of 
mathematics. (Below a connection of exterior differential forms with 
different branches of mathematics will be demonstrated).

\subsection*{Nondegenerate transformations}
One of the fundamental methods in the theory of exterior 
differential forms is application of {\it nondegenerate} 
transformations (below it will be said about {\it degenerate} 
transformations). 

As it has been already noted, in the theory of exterior differential 
forms the nondegenerate transformations are those that conserve the 
differential. This is connected with the property of closed 
differential forms. Since a closed form is a differential (a total one, 
if the form is exact, or an interior one on pseudostructure, if the form 
is inexact), it is evident that the closed form turns out to 
be invariant under all transformations that conserve a differential.

The examples of nondegenerate transformations in the theory  of exterior 
differential forms are unitary, tangent, canonical, gradient transformations.  

To the nondegenerate transformations there are assigned closed forms 
of given degree. To the unitary transformations it is assigned (0-form), to the 
tangent and canonical transformations it is assigned (1-form), to the gradient 
transformations it is assigned (2-form) and so on. 
It should be noted that these transformations are {\it gauge transformations} 
for spinor, scalar, vector, tensor (3-form) fields. 

The above listed transformations are fundamental nondegenerated transformations 
occuring in various branches of mathematics. 

A connection between nondegenerate transformations and closed exterior forms 
disclose an internal commonness of  nondegenerate transformations: all these  
transformations are transformations that preserve a differential. 

One can see that nondegenerate transformations can be classified by a degree of 
coresponding closed differential or dial forms. 

From description of operators of exterior differential forms one can see that 
those are operators that execute some transformations. All these transformations 
are connected with the above listed nondegenerate transformations of exterior 
differential forms.

The possibility to apply nondegenerate transformations shows that the exterior 
differential forms possess  group properties. This extends the 
utilitarian potentialities of the exterior differential forms.

A significance of the nondegenerate transformations consists in the fact that 
they allow one to get new closed differential forms, which 
gives the opportunity to obtain new structures.

\bigskip

Thus, one can see that the skew-symmetric closed exterior differential forms 
possess invariant properties. The invariant properties are a result of 
conjugacy between some objects, to this it points out a closure of exterior 
differential forms. 

With these properties there are connected peculiarities of the mathematical 
apparatus of closed exterior forms, the principal elements of which are 
nondegenerate transformations and identical relations.

These properties of closed exterior differential forms and  peculiarities of 
their mathematical apparatus lie at the basis of practically all branches 
of mathematics. 

Below we outline connections between the exterior differential  
forms and various branches of mathematics in order to show 
what role the exterior differential forms play in mathematics and 
to draw attention to their great potentialities.

\subsection*{Connection between exterior differential forms 
and various branches of mathematics} 

1. {\bf Algebraic and geometrical properties of exterior differential forms}

The basis of Cartan's method of exterior forms, namely, the method 
of analyzing the system of differential equations and manifolds, 
is the basis of the Grassmann algebra (the exterior algebra) [5]. The 
mathematical apparatus of exterior differential forms extends the 
algebraic mathematical apparatus. The differential forms treated 
as elements of algebra allow  studying the manifold structure, 
finding the manifold invariants. Group properties of exterior 
differential forms (a connection with the Lie groups) enable one to 
study the integrability of differential equations. They can make up 
the basis of the invariant field theory. 

The exterior differential forms elucidate an internal connection 
between algebra and geometry. From the definition of the closed 
inexact form it follows that the closed inexact form is a quantity that is 
conservative on pseudostructure. That is, the closed inexact form is 
a conjugacy of algebraic and geometrical approaches. The set 
of relations, namely, the closure condition for exterior form 
($d_\pi \theta^p=0$) and the closure condition for the dual form 
($d_\pi {}^*\theta^p=0$) allow us to describe such a conjugacy. The 
closed form possesses algebraic (invariant) properties, and 
the closed dual form has  geometrical (covariant) properties. 

The mathematical apparatus of exterior differential forms allows 
one to study the elements of interior geometry. This is a fundamental 
formalism in the differential geometry [2]. It allows one to investigate 
the manifold structure. Cartan developed the method of exterior 
forms for investigation of manifolds. It is known that in the case 
of differentiable manifolds the metrical and differential 
characteristics are consistent ones. Such manifolds may be 
regarded as  objects of interior geometry, i.e. there is 
a possibility to study their characteristics as the properties of the surface 
itself irrelevant to imbedding into space. Cartan's structure equation 
is a key tool in studying the manifold structures and fiber spaces [2,3]. 

2. {\bf Theory of functions of complex variables}

The residue method in the theory of analytical functions of complex variables 
is based on the integral theorems by Stokes and 
Cauchy-Poincare that allow us to replace the integral of closed 
form along any closed loop by the integral of this form along 
another closed loop that is homological to the first one [5]. 

As it was already noted, the harmonic function (a differential on 
the Riemann surface) is a closed exterior form:
$$
\theta=pdx+qdy,\quad d\theta=0
$$
to which there corresponds the dual form:
$$
{}^*\theta=-qdx+pdy,\quad d{}^*\theta=0
$$
where $d\theta=0$, $d{}^*\theta=0$ are the Cauchy-Riemann
relations.

3. {\bf Differential equations}

On the basis of the theory of exterior differential forms the 
methods of studying the integrability of the 
system of differential equations, the Pfaff equations (the 
Frobenius theorem), and of finding the integral surfaces have been developed. 
This problem was considered in many works concerning exterior 
differential forms [2, 4, 6].

The operator $d$ appears to be useful for expressing 
the integrability conditions for the systems of partial 
differential equations.

An example of applying the theory of exterior differential forms for 
analyzing integrability of differential equations and determining the 
functional properties of the solutions to these equations is presented in 
Appendix.

4. {\bf Connection with the differential and integral calculus}

As it was already pointed out, the exterior differential forms were 
introduced for designation of integrand expressions that can form integral 
invariants [1].

The exterior differential forms are connected with multiple 
integrals (see, for example, [8]).

It should be mentioned that the closure property of the form 
$f(x)dx$ indicates an existence of antiderivative of the function 
$f(x)$.

Integral relations (1.17), from which it follows the formulas by 
Newton, Liebnitz, Green, 
Gauss-Ostrogradskii, Stokes were derived, are of great importance. Such 
potentials as the Newtonian one, the  potentials of simple and double 
layers are integrals of closed inexact forms. 

The exterior differential forms extend  potentialities of the 
differential and integral calculus.

The theory of integral calculus establishes a connection between 
the differential form calculus and the homology of manifolds 
(cohomology).

5. {\bf Connection with tensors}

As it is well-known, the exterior differential form is a skew-symmetric 
tensor field. 

Information on a connection of the exterior differential 
forms with tensors can be found in [4,8,10,16].

The method of presentation of skew-symmetric tensors as differential forms 
extends possibilities of the mathematical apparatus based on these 
tensors. The tensors are known to be introduced as objects 
that are transformed according to some fixed rule under transforming the 
coordinates. The tensors are attached to the basis that can be transformed 
in an arbitrary way under transition to a new coordinate map. 
The exterior differential form is connected with differentials of 
the coordinates that vary according to the interior characteristics of 
the manifold under translation along manifold. Using differentials of 
coordinates instead of the base vectors enables one directly to make use of the 
integration and differentiation apparatus in physical applications. Instead of 
differentials of coordinates,  a system of linearly-independent exterior 
one-forms can be chosen as the basis, and this makes 
the description independent of the choice of coordinate system [4, 10]. 

As it was pointed out, the operator of exterior 
differentiating $d$ is an abstract generalization of the gradient, 
curl and divergence. This property elucidates a connection 
between the vectorial, algebraic and potential fields. A 
property of the exterior differential form, namely, the existence of 
differential and integral relations between the forms of sequential 
degrees, allows us to classify these fields according to the exterior form 
degree. This was shown for the three-dimensional Euclidean space. 

\bigskip
Thus, even from this brief description of  properties and 
specific features of the exterior differential forms and their 
mathematical apparatus one can clearly see their connection with such 
branches of mathematics as  algebra, geometry, mathematical 
analysis, tensor analysis, differential geometry, differential 
equations, group theory, theory of transformations and so on. 
This is indicative of wide functional and utilitarian potentialities 
of exterior differential forms. One can show that, in essence, the field 
theory is based on invariant and structural properties of exterior 
differential forms that correspond to the conservation laws. 
The exterior differential forms allow us to see an internal connection 
between various branches of mathematics and physics.

The uniqum role of the exterior skew-symmetric differential forms in 
mathematics, as one can see, is connected with their invariant properties. 

Below it will be shown that there exist skew-symmetric differential forms, 
which do not possess the invariant properties, nevertheless their mathematical 
apparatus occurs to be significantly wider. This is due to the fact that these 
differential forms, which were named  evolutionary ones (since, as it will be 
shown below, they possess the evolutionary properties), can generate closed 
exterior differential forms. Such properties of the evolutionary differential 
forms allows a description of discrete transitions, quantum steps, generation 
of various structures and so on,  which cannot be performed within the 
framework of the existing mathematical theories. 

\section{Evolutionary differential forms}

The differential forms, which were named evolutionary ones and differ in their 
properties from the exterior differential forms, may be obtained while studying 
the problem of integrability of differential equations. This can be seen by 
the example of the first-order partial differential equation.
Let
$$ F(x^i,\,u,\,p_i)=0,\quad p_i\,=\,\partial u/\partial x^i \eqno(2.1)$$
be a partial differential equation of the first order. 
Let us consider the functional relation 
$$ du\,=\,\theta\eqno(2.2)$$
where $\theta\,=\,p_i\,dx^i$. 
Here $\theta\,=\,p_i\,dx^i$ is the differential form of the first degree. 
The equation (2.1) will be integrable, if the functional relation is 
an identical one, namely, if the differential form $\theta$ is a  
differential (a closed form). 
In the general case, from equation (2.1) it does not follow (explicitly) 
that the derivatives $p_i\,=\,\partial u/\partial x^i $ that obey 
to the equation (and given boundary or initial conditions of the problem) 
make up a differential, that is, a closed exterior differential form. 
The form $\theta\,=\,p_i\,dx^i$ appears to be an unclosed form and  
is not an identical relation. 
In Appendix it will be shown what is a role of such unclosed forms and  
nonidentical relations in qualitative investigation of functional 
properties of the solutions to differential equations. 

The form $\theta\,=\,p_i\,dx^i$ is an example of the evolutionary differential 
form. (The evolutinary properties of such a form are connected with the 
topological properties of this form commutator, which is nonzero).

Such differential forms originate while investigating the integrability 
of any differential equations that describe various processes 
[12,13].  

\bigskip
The evolutionary differential forms, as well as the exterior forms, are 
skew-symmetric differential forms. 
    
A radical  distinction between the evolutionary forms and the exterior 
ones consists in the fact that the exterior differential forms are defined on 
manifolds with {\it closed metric forms}, whereas the evolutionary differential 
forms are defined on manifolds with {\it unclosed metric forms}. 
 
Before going to description of  evolutionary differential forms it should 
dwell on the properties of manifolds on which skew-symmetric differential 
forms are defined.

\subsection*{Some properties of manifolds}

In the definition of exterior differential forms a differentiable manifold 
was mentioned. Differentiable manifolds 
are topological spaces that locally behave like Euclidean spaces [10,9].

But differentiable manifolds are not a single type of manifolds on which 
the exterior differential forms are defined. In the general case there are 
manifolds with structures of any types. The theory of exterior 
differential forms was developed just for such manifolds. They may 
be the Hausdorff manifolds, fiber spaces, the comological, 
characteristical, configuration manifolds and so on. These manifolds 
and their properties are treated in [2,4,8] and in 
some other works. Since all these manifolds possess 
structures of any types, they have one common property, namely, 
locally they admit one-to-one mapping into the Euclidean subspaces 
and into other manifolds or submanifolds of the same dimension [9]. 

When describing any processes in terms of differential equations, one has to 
deal with manifolds, which do not allow one-to-one 
mapping described above. Such manifolds are, for example, manifolds formed by 
trajectories of elements of the system described by differential equations. 
The manifolds that can be called accompanying manifolds are variable deforming 
manifolds. 
The evolutionary differential forms can be defined on manifolds of this type.

What are the characteristic properties and specific features of 
deforming manifolds and of evolutionary differential forms connected with them? 

To answer this question, let us analyze some properties of metric forms. 

Assume that on the manifold one can set the 
coordinate system with base vectors $\mathbf{e}_\mu$ and define 
the metric forms of manifold [17]: $(\mathbf{e}_\mu\mathbf{e}_\nu)$,
$(\mathbf{e}_\mu dx^\mu)$, $(d\mathbf{e}_\mu)$. The metric forms 
and their commutators define the metric and differential 
characteristics of the manifold.

If metric forms are closed
(the commutators are equal to zero), then the metric is defined 
$g_{\mu\nu}=(\mathbf{e}_\mu\mathbf{e}_\nu)$ and the results of 
translation over manifold of the point 
$d\mathbf{M}=(\mathbf{e}_\mu dx^\mu)$ and of the unit frame 
$d\mathbf{A}=(d\mathbf{e}_\mu)$ prove to be independent of the  
curve shape (the path of integration).

The closed metric forms define the manifold structure, and the commutators of 
metric forms define the manifold differential characteristics that specify 
the manifold deformation: bending, torsion, rotation, twist.

It is evident that the manifolds, that are metric ones or possess the 
structure, have closed metric forms. It is with such manifolds that the 
exterior differential forms are connected.

If the manifolds are deforming manifolds, this means that their 
metric form commutators are nonzero. That is, the metric forms of such 
manifolds turn out to be unclosed. 
The accompanying manifolds and manifolds appearing to be deforming ones are 
examples of such manifolds.

The skew-symmetric evolutionary differential forms 
whose basis are deforming manifolds are defined 
on manifolds with unclosed metric forms.

Thus, the exterior differential forms are skew-symmetric differential forms  
defined on manifolds, submanifolds or on structures with closed 
metric forms. The evolutionary differential forms are skew-symmetric 
differential forms defined on manifolds with metric forms that are unclosed. 

What are the characteristic properties and specific features of such 
manifolds and the related differential forms?

For description of the manifold differential characteristics 
and, correspondingly, the metric forms commutators one can use  
connectednesses [2,4,5,17].

Let us consider the affine connectednesses and their relations to 
commutators of metric forms. 

The components of metric forms can be
expressed in terms of connectedness $\Gamma^\rho_{\mu\nu}$ [17]. The
expressions $\Gamma^\rho_{\mu\nu}$,
$(\Gamma^\rho_{\mu\nu}-\Gamma^\rho_{\nu\mu})$,
$R^\mu_{\nu\rho\sigma}$ are components of the commutators of 
metric forms of zero- first- and third degrees. (The commutator of
the second degree metric form is written down in a more complex
manner [17], and therefore it is not given here).

As it is known [17],
for the Euclidean manifold these commutators vanish identically.
For the Riemann manifold the commutator of the third-degree metric
form is nonzero: $R^\mu_{\nu\rho\sigma}\ne 0$. 

Commutators of metric form vanish in the case of manifolds
that allow local one-to-one mapping into subspaces of
the Euclidean space. In other words, the metric forms of such
manifolds turn out to be closed.

The topological properties of manifolds are connected with the metric 
form commutators. The metric form commutators 
specify a manifold distortion. For example, the commutator of 
the zero degree metric form $\Gamma^\rho_{\mu\nu}$ characterizes the 
bend, that of the first degree form 
$(\Gamma^\rho_{\mu\nu}-\Gamma^\rho_{\nu\mu})$ characterizes the torsion, 
the commutator of the third-degree metric form $R^\mu_{\nu\rho\sigma}$ 
determines the curvature.

The evolutionary properties of differential forms are just connected with 
properties of the metric form commutators. 

\subsection*{Specific features of the evolutionary differential form} 
 
Let us point out some properties of evolutionary forms and show in what their 
difference from  exterior differential forms is manifested. 
The evolutionary differential form of degree $p$ ($p$-form), 
as well as the exterior differential form, can be written down as 
$$
\omega^p=\sum_{\alpha_1\dots\alpha_p}a_{\alpha_1\dots\alpha_p}dx^{\alpha_1}\wedge
dx^{\alpha_2}\wedge\dots \wedge dx^{\alpha_p}\quad 0\leq p\leq n\eqno(2.3)
$$
where the local basis obeys the condition of exterior multiplication 
$$
\begin{array}{l}
dx^{\alpha}\wedge dx^{\alpha}=0\\
dx^{\alpha}\wedge dx^{\beta}=-dx^{\beta}\wedge dx^{\alpha}\quad
\alpha\ne \beta
\end{array}
$$
But the evolutionary form differential cannot be written similarly to that 
presented for exterior differential forms (see formula (1.3)). In the 
evolutionary form differential there appears an additional term connected with 
the fact that the basis of the form changes. For differential forms defined 
on the manifold with unclosed metric form one has 
$d(dx^{\alpha_1}dx^{\alpha_2}\dots dx^{\alpha_p})\neq 0$ (it should be noted 
that for differentiable manifold the following is valid:   
$d(dx^{\alpha_1}dx^{\alpha_2}\dots dx^{\alpha_p}) = 0$).  
For this reason a differential of the evolutionary form $\omega^p$ can be 
written as 
$$
d\omega^p{=}\!\sum_{\alpha_1\dots\alpha_p}\!da_{\alpha_1\dots\alpha_p}dx^{\alpha_1}dx^{\alpha_2}\dots
dx^{\alpha_p}{+}\!\sum_{\alpha_1\dots\alpha_p}\!a_{\alpha_1\dots\alpha_p}d(dx^{\alpha_1}dx^{\alpha_2}\dots
dx^{\alpha_p})\eqno(2.4)
$$
where the second term is connected with a differential of the basis. The second 
term is expressed in terms of the metric form commutator. For the manifold with 
a closed metric form this term vanishes.

For example, we again inspect the first-degree form 
$\omega=a_\alpha dx^\alpha$. The differential of this form can 
be written as $d\omega=K_{\alpha\beta}dx^\alpha dx^\beta$, where 
$K_{\alpha\beta}=a_{\beta;\alpha}-a_{\alpha;\beta}$ are 
components of the commutator of the form $\omega$, and 
$a_{\beta;\alpha}$, $a_{\alpha;\beta}$ are the covariant 
derivatives. If we express the covariant derivatives in terms of 
the connectedness (if it is possible), then they can be written 
as $a_{\beta;\alpha}=\partial a_\beta/\partial
x^\alpha+\Gamma^\sigma_{\beta\alpha}a_\sigma$, where the first 
term results from differentiating the form coefficients, and the 
second term results from differentiating the basis. (In the 
Euclidean space covariant derivatives coincide with ordinary ones 
since in this case derivatives of the basis vanish). If 
we substitute the expressions for covariant derivatives into the 
formula for the commutator components, we obtain the following expression 
for the commutator components of the form $\omega$:
$$
K_{\alpha\beta}=\left(\frac{\partial a_\beta}{\partial
x^\alpha}-\frac{\partial a_\alpha}{\partial
x^\beta}\right)+(\Gamma^\sigma_{\beta\alpha}-
\Gamma^\sigma_{\alpha\beta})a_\sigma\eqno(2.5)
$$
Here the expressions
$(\Gamma^\sigma_{\beta\alpha}-\Gamma^\sigma_{\alpha\beta})$
entered into the second term are just the components of 
commutator of the first-degree metric form.

That is, the corresponding 
metric form commutator will enter into the differential form commutator.

Thus, differentials and, correspondingly, the commutators of exterior and 
evolutionary forms are of different types.

\subsection*{Unclosure of evolutionary differential forms}

The evolutionary differential form commutator, in contrast to that of the 
exterior one, cannot be equal to zero because it involves the metric form 
commutator being nonzero. This means that the evolutionary form differential is 
nonzero. Hence, the evolutionary differential form, in contrast to the case of 
the exterior form,  cannot be closed.  

'he commutators of evolutionary forms depend 
not only on the evolutionary form coefficients but on the characteristics 
of manifolds, on which this form is defined, as well. As a result, such a 
dependence of the evolutionary form commutator produces the topological 
and evolutionary properties of both the commutator and the evolutionary form 
itself (this will be demonstrated below).

Since the evolutionary differential forms are unclosed, the 
mathematical apparatus of evolutionary differential forms does not seem to 
possess any possibilities connected with the algebraic, group, invariant 
and other properties of closed exterior differential forms. However, the 
mathematical apparatus of evolutionary forms includes some new unconventional 
elements. 

Specific features of the mathematical apparatus of evolutionary  
differential forms are shown below. 

\subsection*{Nonidentical relations of evolutionary differential forms}

Above it was shown that the identical relations lie at the 
basis of the mathematical apparatus of exterior differential forms.

In contrast to this, nonidentical relations lie at the basis of the 
mathematical apparatus of evolutionary differential forms. 

The identical relations of closed exterior differential forms reflect 
a conjugacy of any objects. The evolutionary forms, being unclosed, 
cannot directly describe a conjugacy of any objects. But they allow 
description of the process in which the conjugacy may appear (the 
process when closed exterior differential forms are generated). Such 
a process is described by nonidentical relations.

The concept of ``nonidentical relation"  may appear to be inconsistent.
However, it has a deep meaning.

The identical relations establish exact correspondence between the quantities 
(or objects) involved into this relation. It is possible in the case
when the quantities involved into the relation are measurable ones. [A quantity 
is called a measurable quantity if its value does not change under 
transition to another, equivalent, coordinate system. In other words, 
this quantity is invariant one.] In the nonidentical relations one of the 
quantities is unmeasurable. (Nonidentical relations with two unmeasurable 
quantities are meaningless). If this relation is evolutionary one, it turns out 
to be a selfvarying relation, namely, a variation of some object  
leads to  variation of other one, and in turn a variation of 
the second object leads to variation of the first and so on. Since in the 
nonidentical relation one of the objects is a unmeasurable quantity, the 
other object cannot be compared with the first, and therefore, 
the process cannot stop. Here the specific feature is that in the 
process of such selfvarying it may be realized the additional conditions 
under which the identical relation can be obtained from nonidentical relation. 
The additional condition may be realized spontaneously while 
selfvarying the nonidentical relation if the system possesses 
any symmetry. When such additional relations are realized the 
exact correspondence between the quantities involved in the 
relation is established. Under the additional condition a 
unmeasurable quantity becomes a measurable quantity as well, and  
the exact correspondence between the objects 
involved in the relation is established. That is, an identical relation 
can be obtained from a nonidentical relation.

The nonidentical relation is a relation between a closed exterior differential 
form, which is a differential and is a measurable quantity, and 
an evolutionary form, which is an unmeasurable quantity.

Nonidentical relations of such type appear in descriptions of any  
processes (see, for example, functional relation (2.2) built of 
derivatives of differential equation). These relations may be written as 
$$
d\psi \,=\,\omega^p \eqno(2.6)
$$
Here $\omega^p$ is the $p$-degree evolutionary form that is 
nonintegrable, $\psi$ is some form of degree $(p-1)$, and 
the differential $d\psi$ is a closed form of degree $p$.

Nonidentical relations  obtained while describing any processes are 
evolutionary ones.

In the left-hand side of this relation it stands the form differential, 
i.e. a closed form that is an invariant object. In the right-hand 
side it stands the nonintegrable unclosed form that is not an invariant object. 
Such a relation cannot be identical.

One can see a difference of relations for exterior forms and evolutionary ones. 
In the right-hand side of identical relation (see relation (1.19)) 
it stands a closed form, whereas the form in the right-hand side of 
nonidentical relation (2.6) is an unclosed one. 

{\it How is this relation obtained?}
 
Let us consider this by the example of the first degree differential forms. 
A differential of the function of more than one variables can be an example 
of the first degree form. In this 
case the function itself is the exterior form of zero degree.  The state 
function that specifies a state of material system can serve as 
an example of such function. When the physical processes in a material system 
are being described, the state function may 
be unknown, but its derivatives may be known. The values of 
the function derivatives may be equal to some expressions that 
are obtained from the description of a real physical process. And 
it is necessary to find the state function. 

Assume that $\psi$ is the desired state function that depends on the 
variables $x^\alpha$ and also assume that its derivatives in 
various directions are known and equal to some quantities 
$a_\alpha$, namely:
$$
\frac{\partial \psi}{\partial x^\alpha}=a_\alpha\eqno(2.7)
$$
Let us set up the differential expression 
$(\partial\psi/\partial x^\alpha)dx^\alpha$. This differential expression 
is equal to
$$
\frac{\partial\psi}{\partial x^\alpha}dx^\alpha=a_\alpha dx^\alpha\eqno(2.8)
$$
Here the left-hand side of the expression is a differential of 
the function $d\psi$, and in the right-hand 
side it stands the  differential form $\omega=a_\alpha
dx^\alpha$. Relation (2.8) can be written as 
$$
d\psi=\omega\eqno(2.9)
$$
It is evident that relation (2.9) is of the same type as 
(2.6) under the condition that the differential form degrees are equal to 1. 

This relation is nonidentical because the differential form $\omega $ 
is an unclosed differential form. The commutator of this form is nonzero 
since the expressions $a_\alpha $ for the derivatives 
$(\partial\psi/\partial x^\alpha)$ are nonconjugated quantities. 
They are obtained from the description 
of an actual physical process and are unmeasurable quantities.  
(While finding the state function it is commonly assumed that  
its derivatives are conjugated quantities, that is, their mixed 
derivatives are commutative. But in physical processes the 
expressions for these derivatives are usually obtained 
independently of one another. And they appear to be unmeasurable 
quantities, and hence they are not conjugated).

One can come to relation (2.9) by means of analyzing the integrability 
of the partial differential equation. An equation is integrable 
if it can be reduced to the form $d\psi=dU$. However it 
appears that, if the equation is not subjected to an additional
condition (the integrability condition), it is reduced to the 
form (2.9), where $\omega$ is an unclosed form and it cannot be 
expressed as a differential. The first principle of thermodynamics is 
an example of nonidentical relation. 

It arises a question of {\it how to work with nonidentical relation?} 

Two different approaches are possible. 

The first, evident, approach is to find a condition under which the 
nonidentical relation becomes identical and to obtain a closed form under 
this condition. In other words, the nonidentical relation 
is subjected the condition, under which this relation is transformed into 
an identical relation (if it is allowed). 

Such an approach is traditional and is always used implicitly. It may be shown 
that additional conditions are imposed on the mathematical physics equations 
obtained in description of the physical processes so that these equations 
should be invariant (integrable) or should have invariant solutions.

\{Here a psychological point should be noted. While investigating real 
physical processes one often faces the relations that are nonidentical. 
But it is commonly believed that only identical relations can have any 
physical meaning. For this reason one immediately attempts to impose 
a condition onto the nonidentical relation under which this relation 
becomes identical, and it is considered only that this relation can satisfy the 
additional conditions. And all remaining is rejected. 
It is not taken into account that a nonidentical relation is often 
obtained from a description of some physical process and it has  
physical meaning at every stage of the physical process rather 
than at the stage when the additional conditions are 
satisfied. In essence, the physical process does not considered completely. 
At this point it should be emphasized that the nonidentity of 
the evolutionary relation does not mean the imperfect accuracy 
of the mathematical description of a physical process. The 
nonidentical relations are indicative of specific features of 
the physical process development.\} 

This approach does not solve the evolutionary problem. 

Below we present the evolutionary approach to investigation of  
the nonidentical relation. At the basis of this approach 
it lies a peculiarity of the relation that is a nonidentical evolutionary 
relation. Such a relation is a selfvarying relation. 

\subsection*{Selfvariation of the evolutionary nonidentical relation} 
The evolutionary nonidentical relation is selfvarying, 
because, firstly, it is nonidentical, namely, it contains 
two objects one of which appears to be unmeasurable, and, 
secondly, it is an evolutionary relation, namely, a variation of 
any object of the relation in some process leads to variation of 
another object and, in turn, a variation of the latter leads to 
variation of the former. Since one of the objects is an 
unmeasurable quantity, the other cannot be compared 
with the first one, and hence, the process of mutual variation  
cannot stop. 

Varying the evolutionary form coefficients leads to varying the first 
term of the evolutionary form commutator (see (2.5)). In accordance  with this 
variation it varies the second term, that is, the metric form of 
manifold varies. Since the metric form commutators 
specifies the manifold differential characteristics, which are connected with 
the manifold deformation (as it has been pointed out, the commutator of the 
zero degree metric form specifies the bend, that of 
second degree specifies various types of rotation, that of the third 
degree specifies the curvature), this points to the manifold deformation. 
This means that it varies the evolutionary form basis. In turn, this leads to 
variation of the evolutionary form, and the process of intervariation of the 
evolutionary form and the basis is repeated. Processes of variation of 
the evolutionary form and the basis are governed by the evolutionary form 
commutator and it is realized according to the evolutionary relation.

Selfvariation of the evolutionary relation is executed  
by exchange between the evolutionary form coefficients and the manifold 
characteristics. (This may be, for example, a mutual exchange between 
physical quantities and space-time characteristics, or between  
algebraic and geometrical characteristics.)  
This is an exchange between quantities of different nature.
   
The process of the evolutionary relation selfvariation cannot come to the end. 
This is indicated by the fact that both the evolutionary form commutator 
and the evolutionary relation involve unmeasurable quantities.  

The question arises  whether the identical relation can be obtained from this 
nonidentical selfvarying relation, which would allow determination of the 
desired differential form under the differential sign.

It appears that it is possible under the degenerate transformation.

\subsection*{Degenerate transformations.} 
To obtain the identical relation from the evolutionary nonidentical relation,  
it is necessary that a closed exterior differential form should be derived 
from the evolutionary differential form that is included into evolutionary  
relation.   

However, as it was shown above, the evolutionary form cannot be a closed form. 
For this reason a transition from the evolutionary form is possible only 
to an {\it inexact} closed exterior form that is defined on pseudostructure. 

To the pseudostructure there corresponds a closed dual form (whose differential 
vanishes). For this reason a transition from the evolutionary form to a closed 
inexact exterior form proceeds only when the conditions of vanishing the dual 
form differential are realized, in other words, when the metric form 
differential or commutator becomes equal to zero. 

Conditions of vanishing the dual form differential (additional conditions) 
determine the closed metric form and thereby specify the pseudostructure 
(the dual form). In this case the closed exterior ({\it inexact}) form is 
formed. 

Since the evolutionary form differential is nonzero, whereas the closed exterior 
form differential is zero,  a passage from the evolutionary form to the 
closed exterior form is allowed only under {\it degenerate transformation}. 
The conditions of vanishing the dual form differential (an additional 
condition) are the conditions of degenerate transformation. 
 
At this point it should be emphasized that differential, which equals zero, 
is an interior one. The evolutionary form commutator becomes to be zero 
only on the pseudostructure. The total evolutionary form commutator is nonzero. 
That is, under degenerate transformation the evolutionary form differential 
vanishes only {\it on pseudostructure}. The total differential of the 
evolutionary form is nonzero. The evolutionary form remains to be unclosed.   

The conditions of degenerate transformation (additional conditions) can be 
realized, for example, if it will appear any 
symmetries of the evolutionary form coefficients  or its commutator. This can 
happen under selfvariation of the nonidentical relation. (While describing  
material system such additional conditions are related, for example, to 
degrees of freedom of the material system). 

Mathematically to the conditions of degenerate transformation there corresponds 
a requirement that some functional expressions become equal to zero. 
Such functional expressions are Jacobians, determinants, the Poisson brackets, 
residues, and others. 

\subsection*{Obtaining identical relation from  nonidentical one}

Let us consider nonidentical evolutionary relation (2.6).

As it has been already mentioned, the evolutionary differential form $\omega^p$,
involved into this relation is an unclosed one. The commutator, 
and hence the differential, of this form is nonzero. That is, 
$$
d\omega^p\ne 0\eqno(2.10)
$$
If the transformation is degenerate, from the unclosed evolutionary form it 
can be obtained a differential form closed on pseudostructure.  
The differential of this form equals zero. That is, it is 
realized the transition 

 $d\omega^p\ne 0 \to $ (degenerate transformation) $\to d_\pi \omega^p=0$, 
$d_\pi{}^*\omega^p=0$  

On the pseudostructure $\pi$ evolutionary relation (2.6) transforms into 
the relation
$$
d_\pi\psi=\omega_\pi^p\eqno(2.11)
$$
which proves to be the identical relation. Indeed, since the form 
$\omega_\pi^p$ is a closed one, on the pseudostructure it turns 
out to be a differential of some differential form. In other words, 
this form can be written as $\omega_\pi^p=d_\pi\theta$. Relation (2.11) 
is now written as
$$
d_\pi\psi=d_\pi\theta
$$
There are differentials in the left-hand and right-hand sides of
this relation. This means that the relation is an identical one.
                 								 
The evolutionary relation (2.6) becomes 
identical on the pseudostructure. From this relation one can find the desired 
differential $d_\pi\psi$ that is equal to the closed exterior differential 
form derived. One can obtain the desired form $\psi$ from this differential. 

It can be shown that all identical relations of the exterior differential 
form theory are obtained from nonidentical relations (that contain the 
evolutionary forms) by applying degenerate transformations.

Transition from nonidentical relation (2.6) to identical relation (2.11) 
means the following. Firstly, it is from such a relation that 
one can obtain the differential $d_\pi\psi$ and find the desired 
function $\psi_\pi$ (the potential). And, secondly, an emergence 
of the closed (on pseudostructure) inexact exterior form $\omega_\pi^p$  
(right-hand side of relation (2.11)) points to an origination of the  
conservative object. This object is a conservative quantity (the closed 
exterior form  $\omega_\pi^p$) on the pseudostructure (the dual form 
$^*\omega^p$, which defines the pseudostructure). This object is an example 
of the differential and geometrical structure (G-Structure). Such 
objects were already described in section 1 and were named as Bi-Structure. 

This complex is a new conjugated object. Below it will be shown a relation 
between characteristics of these objects and characteristics of the 
evolutionary differential forms.

Thus, the mathematical apparatus of evolutionary differential forms describes 
a process of generation of the closed inexact exterior differential forms, and 
this discloses a process of origination of a new conjugated object.

The evolutionary process of obtaining the identical relation from the 
nonidentical one and obtaing a closed (inexact) exterior form from the 
unclosed evolutionary form describes a process of conjugating any objects.

\subsection*{Transition from nonconjugated operators to conjugated operators} 

In section 1 it has been shown that the condition of the closure of exterior 
differential forms is a result of the conjugacy of any constituents  of the 
exterior or dual forms (the form elements, components of each element, 
exterior and dual forms, exterior forms of 
various degrees, and others). Conversely, if there is 
any closed form, this points to a presence of one or other type of conjugacy.

The conjugacy of any objects is a conjugacy of one or other type 
connected with the exterior differntial forms. To the conjugated objects, 
including operators, there are assigned some closed exterior or dual form.

Since the identical relations of exterior differential forms is a mathematical 
record of the closure conditions of exterior differential forms and, 
correspondingly, of conjugacy of any objects, the process of obtaining 
the identical relation from nonidentical one (selfmodification of the 
nonidentical evolutionary relation and degenerate transformation) is a 
process of conjugating any objects.
 
It can be seen that the process of conjugating the objects is a mutual 
exchange between the quantities of different nature (for example, between 
the algebraic and geometric quantities, between the physical and spatial 
quantities) and vanishing some functional expressions (Jacobians, determinates 
and so on). This follows from the fact that a selfvariation of the 
nonidentical evolutionary relation and a transition from the nonidentical 
evolutionary relation  to identical one develop as a result of mutual 
variations of the evolutionary form coefficients 
(which have the algebraic nature) and the manifold characteristics (which 
have the geometric nature), and a realization of the degenerate transformation 
with obeying additional conditions. 

The evolutionary differential form is an unclosed form, that is, it is the 
form whose differential is not equal to zero. The differential of the 
exterior differential form equals zero. To the closed exterior  
form there correspond conjugated operators, whereas to the evolutionary form 
there correspond nonconjugated operators. A transition from 
the evolutionary form to the closed exterior form is that from  
nonconjugated operators to conjugated ones. This is expressed 
mathematically as a transition from a nonzero differential (the evolutionary 
form differential is nonzero) to a differential that equals zero (the closed 
exterior form differential equals zero). 'his is effected as a transition 
from one coordinate system to another (nonequivalent) coordinate system.  

Since the conjugated objects, to which there correspond the closed exterior 
forms and the identical relations, are obtained from the evolutionary 
differential forms and the nonidentical relations, it is evident that the 
characteristics of conjugated objects are defined by characteristics of the 
evolutionary differential form commutaors, which control the evolutionary 
process. The evolutionary forms enable us to describe the characteristics of 
Bi-Structure originated.

\subsection*{Characteristics of Bi-Structure}
Since the closed exterior differential form, which corresponds to the 
Bi-Structure arisen, was obtained from the nonidentical relation that 
involves the evolutionary form, it is evident that the Bi-Structure 
characteristics must be connected with those of the evolutionary form and of 
the manifold on which this form is defined, with the conditions of degenerate 
transformation and with the values of commutators of the evolutionary form and 
the manifold metric form.

While describing the mechanism of Bi-Structure origination one can see that 
at the instant when the Bi-Structure originates there appear the following 
typical functional expressions and quantities: (1) the condition of degenerate 
transformation, i.e. vanishing of the interior commutator of the metric form; 
(2) vanishing of the interior commutator of the evolutionary form;  
(3) the value of the nonzero total commutator of 
the evolutionary form that involves two terms, namely, 
the first term is composed of the derivatives of the evolutionary form 
coefficients, and the second term is composed of the derivatives of the 
coefficients of the dual form that is connected with the manifold (here 
we deal with a value that the evolutionary form commutator assumes at the 
instant of Bi-Structure origination). 
They determine the following characteristics of the Bi-Structure. The 
conditions of degenerate transformation, as it was said before, determine the 
pseudostructures. The first term of the evolutionary form commutator determines 
the value of the discrete change (the quantum), which the quantity conserved on 
the pseudostructure undergoes at the transition from one pseudostructure to 
another. The second term of the evolutionary form commutator specifies a  
characteristics that fixes the character of the initial manifold deformation,  
which took place before the Bi-Structure arose. (Spin is such an example).

A discrete (quantum) change of a quantity proceeds in the direction 
that is normal (more exactly, transverse) to the pseudostructure. Jumps of the 
derivatives normal to the potential surfaces are examples of such changes.

Bi-Structure may carry a physical meaning. Such bynary objects are the physical  
structures from which the physical fields are formed. This has been shown by 
the  author in the works [12,13]. 

A connection of Bi-Structure with the skew-symmetric differential forms 
allows to introduce a classification of Bi-structure in dependence on 
parameters that specify the skew-symmetric diferential forms and enter into 
nonidentical and identical relation of the skew-symmetric differential forms. 
To determine these parameters one has to consider the problem of integration 
of the nonidentical evolutionary relation. 

\subsection*{Integration of the nonidentical evolutionary relation} 
Under degenerate transformation from the nonidentical evolutionary relation one 
obtains a relation being identical on pseudostructure. 
Since the right-hand side of such a relation can be expressed in terms of 
differential (as well as the left-hand side), one obtains a relation that 
can be integrated, and as a result he obtains a relation with the differential 
forms of less by one degree.

The relation obtained after integration proves to be nonidentical as well. 
 
The resulting nonidentical relation of degree $(p-1)$ (relation that contains 
the forms of the degree $(p-1)$) can be integrated once again if the 
corresponding degenerate transformation has been realized and the identical 
relation has been formed.

By sequential integrating the evolutionary relation of degree $p$ (in the case 
of realization of the corresponding degenerate transformations and forming 
the identical relation), one can get closed (on the pseudostructure) exterior 
forms of degree $k$, where $k$ ranges from $p$ to $0$. 

In this case one can see that under such integration closed (on the 
pseudostructure) exterior forms, which depend on two parameters, are obtained. 
These parameters are the degree of evolutionary form $p$ 
(in the evolutionary relation) and the degree of created closed forms $k$. 

In addition to these parameters, another parameter appears, namely, the 
dimension of space. If the evolutionary relation generates the closed forms 
of degrees $k=p$, $k=p-1$, \dots, $k=0$, to them there correspond the 
pseudostructures of dimensions $(N-k)$, where $N$ is the space dimension. 
\{It is known that to the closed exterior differential forms of degree $k$  
there correspond skew-symmetric tensors of rank $k$ and to corresponding 
dual forms there do the pseudotensors of rank $(N-k)$, where $N$ is 
the space dimensionality. The pseudostructures correspond to such tensors, 
but only on the space formed.\} 

\subsection*{The properties of pseudostructures and closed exterior forms. 
Forming fields and manifolds} 
As mentioned before, the additional conditions, namely, the conditions 
of degenerate transformation, specify the pseudostructure. But at every stage 
of the evolutionary process it is realized only one element of pseudostructure, 
namely, a certain minipseudostructure. The additional conditions determine a 
direction (a derivative of the function that specifies the pseudostructure) 
on which the evolutionary form differential vanishes. (However, in this case 
the total differential of the evolutionary form is nonzero). The closed exterior 
form is formed along this direction. 

While varying the evolutionary variable the minipseudostructures form 
the pseudostructure.

The example of minipseudoctructure is the wave front. The wave front 
is the eikonal surface (the level surface), i.e. the surface with 
a conservative quantity. A direction that specifies the pseudostructure is 
a connection between the evolutionary and spatial variables. It gives the rate 
of changing the spatial variables. Such a rate is a velocity of the wave 
front translation. While its translation the wave front forms the 
pseudostructure.

Manifolds with closed metric forms are formed by pseudostructures. 
They are obtained from manifolds with unclosed metric forms. In this case 
the initial manifold (on which the evolutionary form is defined) and 
the formed manifold with closed metric forms (on which the closed exterior form 
is defined) are different spatial objects. 

It takes place a transition from the initial manifold with unclosed metric 
form to the pseudostructure, namely, to the created manifold with closed 
metric forms. 
Mathematically this transition (degenerate transformation) proceeds as {\it a 
transition from one frame of reference to another, nonequivalent, frame of 
reference.}

The pseudostructures, on which the closed {\it inexact} forms are defined, 
form the pseudomanifolds. (Integral surfaces, pseudo-Riemann and 
pseudo-Euclidean spaces are the examples of such manifolds). In this process 
dimensions of formed manifolds are connected with the evolutionary 
form degree.

To transition from pseudomanifolds to metric manifolds there corresponds 
a transition from closed {\it inexact} differential forms to {\it exact} 
exterior differential forms. (Euclidean and Riemann spaces are examples 
of metric manifolds).

Here it is to be noted that the examples of pseudometric spaces are  
potential surfaces (surfaces of a simple layer, a double layer and so on). 
In these cases the type of potential surfaces is connected with the above 
listed parameters.

Since the closed metric form is dual with respect to some closed exterior 
differential form, the metric forms cannot become closed by themselves, 
independently of the exterior differential form. This proves that manifolds 
with closed metric forms are connected with the closed exterior differential 
forms. This indicates that the fields of conservative quantities are formed 
from closed exterior forms at the same time when the manifolds are created from 
the pseudoctructures. (The specific feature of the manifolds with closed metric 
forms that have been formed is that they can carry some information.) 
That is, the closed exterior differential forms and manifolds, on which they 
are defined, are mutually connected objects. On the one hand, this shows 
duality of these two objects (the pseudostructure and the closed inexact 
exterior form), and, on the  other hand, this means that these objects 
constitute a unified whole. This whole is a new conjugated object (Bi-Structure). 

\bigskip
{\large\bf Summary}.

In section 1 it had been described invariant properties of skew-symmetric 
closed exterior  differential forms and was shown that due to their properties 
the closed exterior forms play a fundamental role in various branches of 
mathematics. 

In section 2 we showed that besides the exterior skew-symmetric differential 
forms, which possess invariant properties, there are skew-symmetric 
differential forms, which possess the evolutionary properties (which, for 
this reason, have been named the evolutionary forms). 

By comparing the exterior and evolutionary forms one can see that they possess 
the opposite properties. At the same time the exterior and evolutionary forms 
constitute a unified whole: evolutionary differential forms generate the closed 
exterior differential forms. 

The mathematical apparatus of exterior and evolutionary skew-symmetric 
differential forms constitute a new closed mathematical apparatus that 
possesses the unique properties. It includes new, unconventional, elements: 
"nonidentical relation", "degenerate transformation", "transition from one 
frame of reference to another, nonequivalent, frame of reference". This 
allows to create the mathematical language that has radically new abilities.

Due to their properties, the skew-symmetric differential forms enable one 
to see the internal connection between various branches of mathematics. 
Many foundations of the mathematical apparatus of skew-symmetric differential 
forms presented in this work may turn out to be of great importance for various 
sections of mathematics. Identical and nonidentical relations, 
nondegenerate  and degenerate transformations, transitions from nonidentical 
relations to identical ones, transition from nonconjugated operators to 
conjugated operators, a possibility to describe  an origination of structure and 
formation of fields and manifolds, and other foundations and potentialities of 
the mathematical apparatus of skew-symmetric differential 
forms presented may find many applications in such branches of mathematics as 
the qualitative theory of differential and integral equations, differential 
geometry and topology, theory of functions, theory of series, theory of 
numbers and others.  

Due to their properties skew-symmetric differential forms have a great 
significance in applcations. In the works by the author [12,13] it has 
been shown that the apparatus of skew-symmetrical differential forms 
discloses a mechanism of originating the physical structures and forming 
the physical fields, and also explains the causality of these processes.

Below we present the example of application of the skew-symmetric differential 
forms for qualitive studying solutions of differential equations. 

\bigskip
\rightline{\large\bf Appendix }

{\bf Application of the mathematical apparatus of the skew-symmetric 
differential forms to qualitative investigation of functional properties of the 
solutions to differential equations}

\bigskip
The presented method of investigating the solutions to differential equations 
is not new. Such an approach was developed by Cartan [2] in his analysis of the
integrability of differential equations. Here this apptoach is presented to 
demonstrate a role of exterior and evolutionary differential forms.

A role of exterior differential forms in the qualitative
investigation of the solutions to differential equations is conditioned by
the fact that the mathematical apparatus of these forms enables one to
determine the conditions of consistency for various elements of differential
equations or for the system of differential equations. This enables one, for
example, to define the consistence of the partial derivatives in the partial
differential equations, the consistency of the differential equations in the 
system of differential equations, the conjugacy of the function derivatives and
of the initial data derivatives in ordinary differential equations and 
so on. The functional properties of the solutions to differential equations 
are just depend on whether or not the conjugacy conditions are satisfied.

The basic idea of the qualitative investigation of the solutions to
differential equations can be clarified by the example of the first-order
partial differential equation.

Let
$$ F(x^i,\,u,\,p_i)=0,\quad p_i\,=\,\partial u/\partial x^i \eqno(A.1)$$
be the partial differential equation of the first order. Let us consider the
functional relation
$$ du\,=\,\theta\eqno(A.2)$$
where $\theta\,=\,p_i\,dx^i$ (the summation over repeated indices is implied).
Here $\theta\,=\,p_i\,dx^i$ is the differential form of the first degree.

The specific feature of functional relation (A.2) is that in
the general case this relation turns out to be nonidentical.

The left-hand side of this relation involves a differential, and
the right-hand side includes the differential form  $\theta\,=\,p_i\,dx^i$. 
For this relation to be identical, the differential form 
$\theta\,=\,p_i\,dx^i$ must be a differential as well
(like the left-hand side of relation (A.2)), that is, it has to be a
closed exterior differential form. To do this it requires the commutator
$K_{ij}=\partial p_j/\partial x^i-\partial p_i/\partial x^j$ of the
differential form $\theta $ has to vanish.

In the general case, from equation (A.1) it does not follow (explicitly)
that the derivatives $p_i\,=\,\partial u/\partial x^i $ that obey
to the equation (and given boundary or initial conditions of the problem)
make up a differential. Without any supplementary conditions
the commutator of the differential form $\theta $ defined as
$K_{ij}=\partial p_j/\partial x^i-\partial p_i/\partial x^j$ is
not equal to zero. The form $\theta\,=\,p_i\,dx^i$ proves to be unclosed and 
is not a differential like the left-hand side of relation (A.2). The functional
relation (A.2) appears to be nonidentical: the left-hand side of this relation
is a differential, but the right-hand side is not a differential. 

Functional relation (A.2) is an example of nonidentical evolutionary relation.  

The nonidentity of functional relation (A.2) points to a fact
that without additional conditions derivatives of the initial
equation do not make up a differential. This means that the corresponding
solution to the differential equation $u$ will not be a function of $x^i$.
It will depend on the commutator of the form $\theta $, that is, it will be
a functional.

To obtain the solution that is the function (i.e., derivatives of this
solution form a differential), it is necessary to add the closure condition
for the form $\theta\,=\,p_idx^i$ and for the dual form (in the present
case the functional $F$ plays a role of the form dual to $\theta $) [2]:
$$\cases {dF(x^i,\,u,\,p_i)\,=\,0\cr
d(p_i\,dx^i)\,=\,0\cr}\eqno(A.3)$$
If we expand the differentials, we get a set of homogeneous equations
with respect to $dx^i$ and $dp_i$ (in the $2n$-dimensional space -- initial
and tangential):
$$\cases {\displaystyle \left ({{\partial F}\over {\partial x^i}}\,+\,
{{\partial F}\over {\partial u}}\,p_i\right )\,dx^i\,+\,
{{\partial F}\over {\partial p_i}}\,dp_i \,=\,0\cr
dp_i\,dx^i\,-\,dx^i\,dp_i\,=\,0\cr} \eqno(A.4)$$
The solvability conditions for this set (vanishing of the determinant
composed of coefficients at $dx^i$, $dp_i$) have the form:
$$
{{dx^i}\over {\partial F/\partial p_i}}\,=\,{{-dp_i}\over {\partial F/\partial x^i+p_i\partial F/\partial u}} \eqno (A.5)
$$
These conditions determine an integrating direction, namely, a pseudostructure,
on which the form $\theta \,=\,p_i\,dx^i$ turns out to be closed one,
i.e. it becomes a differential, and from relation (A.2) the identical relation
is produced.
If conditions (A.5), that may be called the integrability conditions,
are satisfied, the derivatives constitute
a differential $\delta u\,=\,p_idx^i\,=\,du$ (on the pseudostructure), and the
solution becomes a function.
Just such solutions, namely, functions on the pseudostructures
formed by the integrating directions, are the so-called generalized
solutions [19]. The derivatives of the generalized solution constitute the
exterior form that is closed on the pseudostructure.

(If conditions (A.5) are not satisfied, that is, the derivatives do not
form a differential, the solution that corresponds to such derivatives
will depend on the differential form commutator constructed of derivatives.
That means that the solution is a functional rather then a function.)

Since the functions that are the generalized solutions
are defined only on the pseudostructures, they have discontinuities in
derivatives in the directions that are transverse to the pseudostructures.
The order of derivatives with discontinuities
is equal to  the exterior form degree. If the form of zero
degree is involved in the functional relation, the function itself,
being a generalized solution, will have discontinuities.

If we find the characteristics of equation (A.1), it appears that
conditions (A.5) are the equations for  characteristics [11].
That is, the characteristics are examples of
the pseudostructures on which  derivatives of the differential equation
constitute the closed forms and the solutions prove to be the functions 
(generalized solutions). (The characteristic manifolds of equation (A.1) are 
the pseudostructures $\pi$ on which the form $\theta =p_idx^i$ becomes a closed
form: $\theta _{\pi}=d u_{\pi}$).

Here it is worth noting that coordinates of the equations for
characteristics are not identical to independent coordinates of 
initial space on which equation (A.1) is defined. A transition from
the initial space to the characteristic manifold appears to be a
{\it degenerate} transformation, namely,
the determinant of the set of equations (A.4) becomes zero. The derivatives
of  equation (A.1) are transformed from the tangent space to the cotangent one.
The transition from the tangent space, where the commutator of the form 
$\theta$ is nonzero (the form is unclosed, the derivatives do not form a 
differential), to the characteristic manifold, namely, the cotangent space, 
where the commutator becomes equal to zero (the closed exterior form is formed, 
i.e. the derivatives form a differential), is the example of the degenerate 
transformation. 

A partial differential equation of the first order has been analyzed,
and the functional relation with the form of the first degree analogous
to the evolutionary form has been considered.

Similar functional properties have the solutions to all differential equations. 
And, if the order of the differential equation is $k$, the functional relation
with the $k$-degree form corresponds to this equation. For ordinary differential 
equations the commutator is produced at the expense of the conjugacy of  
derivatives of the functions desired and those of the initial data (the 
dependence of the solution on the initial data is governed by the commutator).

In a similar manner one can also investigate  the solutions to a set of 
partial differential equations and the solutions to ordinary differential
equations (for which the nonconjugacy of desired functions and initial
conditions is examined).

It can be shown that  the solutions to equations of mathematical
physics, on which no additional external conditions are imposed, are
functionals. The solutions prove to be exact only under realization of
the additional requirements, namely, the conditions of degenerate transformations: 
vanishing determinants, Jacobians and so on, that define the integral
surfaces.  The characteristic manifolds, the envelopes of characteristics,
singular points, potentials of simple and double layers, residues and others 
are the examples of such surfaces.

Here the mention should be made of the generalized Cauchy problem when the
initial conditions are given on some surface. The so called ``unique" solution
to the Cauchy problem, when the output derivatives can be determined (that is,
when the determinant built of the expressions at these derivatives is
nonzero), is a functional since the derivatives obtained in such a way
prove to be nonconjugated, that is, their mixed derivatives form a
commutator with nonzero value, and the solution depends on this commutator.

The dependence of the solution on the commutator may lead to instability of the 
solution. Equations that do not possess the integrability conditions (the 
conditions such as, for example, the characteristics, singular points, 
integrating factors and others) may have the unstable solutions. Unstable 
solutions appear in the case when the additional conditions are not realized 
and no exact solutions (their derivatives form a differential) are formed. 
Thus, the solutions to the equations of the elliptic type may be unstable.

Investigation of nonidentical functional relations lies at the basis
of the qualitative theory of differential equations. It is well known that
the qualitative theory of differential equations is based on the analysis of
unstable solutions and integrability conditions. From the functional
relation it follows that the dependence of the solution on the commutator
leads to instability, and the closure conditions of the forms constructed by
derivatives are the integrability conditions. One can see that
the problem of unstable solutions and integrability conditions appears,
in fact, to be reduced to the question  under what conditions the identical
relation for the closed form is produced from the nonidentical relation that
corresponds to the relevant differential equation (the relation such as (A.2)), 
the identical relation for the closed form is produced. In other words, whether 
or not the solutions are functionals? This is to the same question that the 
analysis of the correctness of setting the problems of mathematical physics is 
reduced.

Here the following  should be emphasized. When the degenerate
transformation from the initial nonidentical functional relation is fulfilled,
an integrable identical relation is obtained. As a result of integrating, 
one obtains a relation that contains exterior forms of less by one
degree and which once again proves to be (in the general case without 
additional conditions) nonidentical. By integrating the functional relations 
sequentially obtained (it is possible only under realization of the degenerate 
transformations) from the initial functional relation of degree $k$ one can obtain 
$(k+1)$ functional relations each involving exterior forms of one of degrees: 
$k, \,k-1, \,...0$. In particular, for the first-order partial differential 
equation it is also necessary to analyze the functional relation of zero degree. 

Thus, application of the exterior differential forms allows one to reveal 
the functional properties of the solutions to differential equations.

It is evident that, for solutions to the differential equation be generalized 
solutions (i.e. solutions whose derivatives form a differential, namely, 
the closed form), the differential equation has to be subject the additional 
conditions. Clearly, only generalized solutions  can correspond to various 
structures including physical structures.
Let us consider what equations are obtained in this case. 

Return to equation (A.1).

Assume that it does not explicitly depend on $u$ and is solved with respect 
to some variable, for example $t$, that is, it has the form of 
$${{\partial u}\over {\partial t}}\,+\,E(t,\,x^j,\,p_j)\,=\,0, \quad p_j\,=\,{{\partial u}\over {\partial x^j}}\eqno(A.6)
$$
Then integrability conditions (A.5)
(the closure conditions of the differential form
$\theta =p_idx^i$  and the corresponding dual form) can be written as
(in this case $\partial F/\partial p_1=1$)
$${{dx^j}\over {dt}}\,=\,{{\partial E}\over {\partial p_j}}, \quad
{{dp_j}\over {dt}}\,=\,-{{\partial E}\over {\partial x^j}}\eqno(A.7)$$

These are the characteristic relations for equation (A.6). As it is well known, 
the canonical relations have just such a form. 

As a result we conclude that the canonical relations
are the characteristics of equation (A.6) and the integrability conditions
for this equation. 

The canonical relations that are obtained from the closure condition of the
differential form $\theta = p_idx^i$ and the corresponding dual form, are
the examples of the identical relation of the theory of exterior differential 
forms.

The conjugated coordinates $x^j$, $p_j$ that obey the
canonical relations can be put into correspondence with another equation,
namely, the equation for the canonical relation integral
$V(t,\,x^j,\,p_j)$
$$ {{\partial V}\over {\partial t}}\,+\,[E,\,V]\,=\,0\eqno (A.7`)$$
where $[E,\,V]$ is the Poisson bracket.

Equation (A.6) provided with the supplementary conditions,
namely, the canonical relations (A.7), is called the Hamilton-Jacobi equation
[11]. In other words, the equation whose derivatives obey the canonical
relation is referred to as the Hamilton-Jacobi equation. The derivatives of
this equation form the differential, i.e. the closed exterior differential form:
$\delta u\,=\,(\partial u/\partial t)\,dt+p_j\,dx^j\,=\,-E\,dt+p_j\,dx^j\,=\,du$.

The equations of field theory belong to this type.
$${{\partial s}\over {\partial t}}+H \left(t,\,q_j,\,{{\partial s}\over {\partial q_j}}
\right )\,=\,0,\quad
{{\partial s}\over {\partial q_j}}\,=\,p_j \eqno(A.8)$$
where $s$ is the field function for the action functional $S\,=\,\int\,L\,dt$.
Here $L$ is the Lagrange function, $H$ is the Hamilton function:
$H(t,\,q_j,\,p_j)\,=\,p_j\,\dot q_j-L$, $p_j\,=\,\partial L/\partial \dot q_j$.
The closed form $ds\,=\,H\,dt\,+\,p_j\,dq_j$ (the
Poincare invariant) corresponds to equation (A.8).

The coordinates $q_j$, $p_j$ in equation (A.8) are the conjugated ones. 
They obey the canonical relations. The equation for the canonical relation 
integral that is similar to equation (A.7`) can be assigned to equation (A.8).

In quantum mechanics (where to the coordinates $q_j$, $p_j$ the operators are
assigned) the Schr\H{o}dinger equation [18] serves as an analog to equation 
(A.8), and the Heisenberg equation serves as an analog to the relevant equation 
for the canonical relation integral. Whereas the closed
exterior differential form of zero degree (the analog to the Poincare
invariant) corresponds to the Schr\H{o}dinger equation, the closed dual form
corresponds to the Heisenberg equation.

A peculiarity of the degenerate transformation can be considered by the example
of the field equation. In section 2 it was said that a transition from the
unclosed differential form (which is included into the functional relation) to 
the closed form is the degenerate transformation. Under degenerate 
transformation a transition from the initial manifold (on which the differential 
equation is defined) to the characteristic (integral) manifold goes on.

Here the degenerate transformation is a transition from the Lagrange
function to the Hamilton function. The equation for the Lagrange function, that
is the Euler variational equation, was obtained from the condition 
$\delta S\,=\,0$, where $S$ is the action functional. In the real case, 
when forces are nonpotential or couplings are nonholonomic, the quantity 
$\delta S$ is not a closed form, that is, $d\,\delta S\,\neq \,0$.
But the Hamilton function is obtained from the condition $d\,\delta S\,=\,0$
which is the closure condition for the form $\delta S$. A transition from the
Lagrange function $L$ to the Hamilton function $H$ (a transition from
variables $q_j,\,\dot q_j$ to variables $q_j,\,p_j=\partial L/\partial \dot q_j$)
is a transition from the tangent space, where the form is unclosed, to
the cotangent space with a closed form. One can see that this transition is
a degenerate one.

The invariant field theories used only nondegenerate transformations that
conserve the differential. By the example of the canonical relations it is 
possible to show that nondegenerate and degenerate transformations 
are connected. The canonical relations in the invariant field
theory correspond to nondegenerate tangent transformations.
At the same time, the canonical relations coincide  with the characteristic
relation for  equation (A.8), which
the degenerate transformations correspond to. The degenerate transformation is
a transition from the tangent space ($q_j,\,\dot q_j)$) to the
cotangent (characteristic) manifold ($q_j,\,p_j$). 
(This is a transition from the manifold that corresponds
to the material system, to physical fields. Such a transition is connected with
emergence of a physical structure). On the other hand, the nondegenerate
transformation is a transition
from one characteristic manifold ($q_j,\,p_j$) to the other characteristic
manifold ($Q_j,\,P_j$). (Physically, this describes a transition from one
physical structure to another physical
structure.) $\{$The formula of canonical transformation can be written as
$p_jdq_j=P_jdQ_j+dW$, where $W$ is the generating function$\}$.

It may be easily shown that such a property of duality is also a specific
feature of transformations such as tangent, gradient, contact, gauge, conform 
mapping, and others.

1. Cartan E., Lecons sur les Invariants Integraux. -Paris, Hermann, 1922.  

2. Cartan E., Les Systemes Differentials Exterieus ef Leurs Application 
Geometriques. -Paris, Hermann, 1945.  

3. Sternberg S., Lectures on Differential Geometry. -Englewood Cliffd , N.J.:
Prentiice-Hall,1964

4. Schutz B.~F., Geometrical Methods of Mathematical Physics. Cambrige 
University Press, Cambrige, 1982.

5. Encyclopedia of Mathematics. -Moscow, Sov.~Encyc., 1979 (in Russian).

6. Finikov S.~P.,  Method of the Exterior Differential Forms by Cartan in 
the Differential Geometry. Moscow-Leningrad, 1948 (in Russian).

7. Bott R., Tu L.~W., Differential Forms in Algebraic Topology. 
Springer, NY, 1982.

8. Novikov S.~P., Fomenko A.~P., Elements of the differential geometry and 
topology. -Moscow, Nauka, 1987 (in Russian). 

9. Efimov N.~V. Exterior Differntial Forms in the Euclidian Space. 
Moscow State Univ., 1971.

10. Wheeler J.~A., Neutrino, Gravitation and Geometry. Bologna, 1960.

11. Smirnov V.~I., A course of higher mathematics. -Moscow, 
Tech.~Theor.~Lit. 1957, V.~4 (in Russian).

12. Petrova L.~I., Exterior and evolutionary skew-symmetric differential forms 
and their role in mathematical physics. 

http://arXiv.org/pdf/math-ph/0310050

13. Petrova L.~I., Conservation laws. Their role in evolutionary processes.   
(The method of skew-symmetric differential forms). 

http://arXiv.org/pdf/math-ph/0311008

14. Haywood R.~W., Equilibrium Thermodynamics. Wiley Inc. 1980.

15. Fock V.~A., Theory of space, time, and gravitation. -Moscow, 
Tech.~Theor.~Lit., 1955 (in Russian).

16. Konopleva N.~P. and Popov V.~N., The gauge fields. Moscow, Atomizdat, 1980 
(in Russian).

17. Tonnelat M.-A., Les principles de la theorie electromagnetique 
et la relativite. Masson, Paris, 1959.

18. Dirac P.~A.~M., The Principles of Quantum Mechanics. Clarendon Press, 
Oxford, UK, 1958.

19. Vladimirov V.~S., Equations of the mathematical physics. -Moscow, 
Nauka, 1988 (in Russian).

\end{document}